\documentclass[11pt,reqno,oneside]{article}
\usepackage[left=1in,right=1in,top=1in,bottom=1in]{geometry}
\usepackage{amsmath,amssymb}
\usepackage{theorem}
\usepackage{color}

\theorembodyfont{\itshape}
\theoremstyle{plain}
\newtheorem{theorem}{Theorem}[section]

\newtheorem{lemma}[theorem]{Lemma}
\newtheorem{fact}[theorem]{Fact}
\newtheorem{corollary}[theorem]{Corollary}
\newtheorem{proposition}[theorem]{Proposition}
\theorembodyfont{\rmfamily}
\newtheorem{example}[theorem]{Example}
\newtheorem{remark}[theorem]{Remark}
\newtheorem{question}[theorem]{Question}
\newtheorem{problem}[theorem]{Problem}
\newtheorem{definition}[theorem]{Definition}

\newtheorem{claim}{Claim}

\newcommand{\id}{\mbox{{\rm id}}}

\def\Z{\mathbb{Z}}
\def\Q{\mathbb{Q}}
\def\R{\mathbb{R}}

\def\inv#1{\Phi({#1})}
\def\invbrackets#1{\Phi\left({#1}\right)}
\def\grp#1{\langle{#1}\rangle}
\def\Aut{\mathrm{Aut}}

\begin{document}

\title{Automorphism groups  of  dense subgroups of $\mathbb R^n$}
\author{Vitalij A.~Chatyrko and Dmitri B.~Shakhmatov}

\date{}

\maketitle
\begin{abstract}
By an {\em automorphism\/} of a topological group $G$
we mean 
an isomorphism of $G$ onto itself which is also a homeomorphism. 
In this article, we 
study the automorphism group $\mathrm{Aut}(G)$ 
of a dense subgroup $G$ of $\mathbb R^n, n \geq 1.$ 
We
show that $\mathrm{Aut}(G)$ can be naturally identified with the subgroup $\inv{G}=\{A\in GL(n,\R): G\cdot A =G\}$ of the group $GL(n,\R)$ of all non-degenerated 
$(n \times n)$-matrices over $\mathbb R$, where 
$G\cdot A=\{g\cdot A:g\in G\}$.
We 
describe $\inv{G}$ for many dense subgroups $G$ of either $\R$ or $\R^2$. We consider also an inverse problem of which symmetric subgroups of $GL(n,\R)$ can be realized as $\inv{G}$ for some dense subgroup $G$ of $\R^n$. For example, for $n\ge 2$, we show that the group $SL^{\pm}(n,\R)=\{A\in GL(n,\R): \mathrm{det} A=\pm 1\}$ cannot be realized in this way. The realization problem is quite non-trivial even in the one-dimensional case
and has deep connections to number theory.
\end{abstract} 


\smallskip
{\it 2010 AMS (MOS) Subj. Class.:} Primary: 22A05; Secondary: 11R04, 12F20, 16S34, 20D45, 20E36, 20K30, 46A16, 47D03, 54F45, 54H11, 54H13.
\bigskip

{\em All (topological) vector spaces considered in this paper are assumed to be over the field $\R$ of real numbers.\/}   

As usual, $\mathbb{N}$ denotes the set of natural numbers, $\Z$ denotes the ring of integer numbers, $\Q$ denotes the field of rational numbers and $\R$ denotes the field of real numbers. The symbol $\R^n$ denotes the $n$-dimensional Euclidean space considered with its norm topology.

For topological groups $G$ and $H$, the notation $G\cong H$ means that $G$ and  $H$ are isomorphic as topological groups.

 For necessary information on topological groups we refer to \cite{AT}.

\section{Introduction}

The symbol $GL(V)$  denotes the group of linear bijective maps of a vector space $V$ 
with the composition of maps as its group operation. When $V$ is a topological vector space, we use the symbol $HGL(V)$ to denote the subgroup of $GL(V)$ consisting of continuous maps $\varphi\in GL(V)$ whose inverse $\varphi^{-1}$ is also continuous; that is,
$HGL(V)$ is the group of linear bijections of $V$ which are simultaneously homeomorphisms of  $V$ onto itself.

Let $X$ be a subspace of a topological vector space $V$. We denote by 
\begin{equation}
\label{eq:1}
\inv{X}=\{\varphi\in HGL(V): \varphi(X)=X\}
\end{equation}
the subgroup of $HGL(V)$ consisting of all $\varphi$ that map $X$ onto itself, or equivalently, of all maps $\varphi\in HGL(V)$ whose restriction $\varphi\restriction_X:X\to X$ to $X$ is a homeomorphism of $X$ onto itself. The restriction map 
\begin{equation}
\label{eq:0}
r_X: \inv{X}\to \mathrm{Homeo}(X)
\end{equation}
  defined by
\begin{equation}
\label{eq:2}
r_X(\varphi)=\varphi\restriction_X
\ 
\text{ for }
\ 
\varphi\in \inv{X}
\end{equation}
is a group homomorphism from $\inv{X}$ to the group $\mathrm{Homeo}(X)$ of homeomorphisms of $X$ onto itself.
To check that $r_X$ is a homomorphism, fix $\varphi,\psi\in \inv{X}$ and note that 
$$
r_X(\varphi \circ \psi^{-1})=(\varphi\circ \psi^{-1})\restriction_X
=
\varphi\restriction_X\circ \psi^{-1}\restriction_X
=
\varphi\restriction_X\circ (\psi\restriction_X)^{-1}
=
r_X(\varphi)\circ r_X(\psi)^{-1}.
$$

\begin{definition}
Let $G$ be a topological group. 
\begin{itemize}
\item[(i)]
A map $f:G\to G$ will be called an {\em automorphism\/} of $G$ provided that $f$ is simultaneously 
an isomorphism (of the group $G$ onto itself) and a homeomorphism (of the topological space $G$ onto itself).
\item[(ii)] The symbol $\mathrm{Aut(G)}$ denotes the group of all automorphisms of $G$ with the composition of maps as its group operation.
\end{itemize}
\end{definition}

\begin{proposition}
\label{basic:proposition}
If $X$ is a subgroup of the additive group $(V,+)$ of a topological vector space $V$, then the map $r_X:\inv{X}\to \Aut(X)$ defined in \eqref{eq:2}
is a group homomorphism from the subgroup $\inv{X}$ of the group $HGL(V)$ to the automorphism group $\Aut(X)$ of $X$.
\end{proposition} 
Proof.
Suppose that $X$ is a subgroup of the additive group $(V,+)$ of $V$.
Since each $\varphi\in\inv{X}$ is a linear bijection of $V$ onto itself and $\varphi(X)=X$ by \eqref{eq:1}, it easily follows that 
$\varphi\restriction_X$ is a bijective homomorphism of $X$ onto itself. 
Since $r_X(\varphi)=\varphi\restriction_X$ is a homeomorphism of $X$ onto itself by 
\eqref{eq:0} and \eqref{eq:2}, we conclude that 
$r_X(\varphi)$ is an automorphism of the topological group $X$, i.e.,
$r_X(\varphi)\in \Aut(X)$. 
$\Box$

\medskip

We shall say that a topological vector space $V$ is {\em complete\/} if its additive group $(V,+)$ is Raikov complete (cf. \cite{AT}).

Quite conveniently, when the subgroup $X$ of $V$ in Proposition \ref{basic:proposition} is dense in a complete vector space $V$, the homomorphism $r_X:\inv{X}\to \Aut(X)$ becomes an isomorphism.
\begin{theorem}
\label{isomorphism:theorem}
Let $G$ be a dense subgroup (of the additive group $(V,+)$) of a complete vector space $V$. Then the map $r_G:\inv{G}\to \Aut(G)$ is an isomorphism between
the subgroup $\inv{G}$ of the group $HGL(V)$ and the automorphism group $\Aut(G)$ of $G$.
In particular, $\inv{G}$ is an isomorphic copy of $\Aut(G)$ in $HGL(V)$.
\end{theorem}

The proof of this theorem will be carried out in Section~\ref{sec:2}.

In this paper we are concerned with two general problems inspired by Theorem~\ref{isomorphism:theorem}.

\begin{problem}
\label{pro:1}
Given a dense subgroup 
$G$
(of the additive group $(V,+)$) of a complete vector space $V$, 
``calculate''
the automorphism group $\Aut(G)$ of $G$, or equivalently, the subgroup $\inv{G}$ of $HGL(V)$.
\end{problem}

The second problem is a converse problem to 
Problem~\ref{pro:1}.
\begin{problem}
\label{pro:2}
Given a complete vector space $V$ and a subgroup $H$ of 
the group $HGL(V)$, can one find a dense subgroup $G$ (of the additive group $(V,+)$) of $V$  such that $\inv{G}=H$? In other words, which subgroups of $HGL(V)$ are naturally isomorphic to the automorphism group $\Aut(G)$ of some dense subgroup $G$ of $V$? 
\end{problem}

In this paper, we make some progress on these two problems in the special case when the complete vector space $V$ in question is finite-dimensional (and thus, isomorphic to $\R^n$ for some positive integer $n$).

When $V=\R^n$ for some positive integer $n$, all linear maps from $V$ to itself are automatically continuous, so the equality $HGL(V)=GL(V)$ holds. Furthermore, in this special case $GL(V)$
is isomorphic to the general linear group $GL(n,\R)$ of non-degenerated $(n\times n)$-matrices with real coefficients, so $\inv{G}$ can be viewed as a subgroup of the matrix group $GL(n,\R)$, and the inverse map $r_G^{-1}: \Aut(G)\to \inv{G}\subseteq GL(V)\cong GL(n,\R)$ of $r_G$ can be viewed as a faithful representation of the automorphism group $\Aut(G)$ in the matrix group $GL(n,\R)$.

We note that even the one-dimensional case $V=\R$ is rather complicated, with Problems~\ref{pro:1} and~\ref{pro:2} having interesting and non-trivial connections with number theory. 

The paper is organized as follows. The proof of Theorem~\ref{isomorphism:theorem} is given in Section~\ref{sec:2}.
Basic properties of the subgroup $\inv{G}$ of $HGL(V)$ are studied in Section~\ref{sec:3}. In particular,
it is proved in Corollary~\ref{elementary:properties:of:Phi(G):corollary}(i) that $\inv{G}$ is symmetric with respect to addition, i.e. satisfies the equation $-\inv{G}=\inv{G}$. 

The case $V=\R$ is studied in Sections~\ref{sec:4}--\ref{sec:6}.
In this case the group $HGL(\R)$ coincides with the multiplicative group $\R^\times$ of non-zero real numbers (Remark~\ref{linear:maps:in:dimension:1}). Another specific fact in the one-dimensional case is that a subgroup of $\R$ is dense in $\R$ if and only if it is non-cyclic (Remark~\ref{non-cyclic:remark}). Moreover, a cyclic subgroup $C$ of $\R$ is rigid; that is, its automorphism group $\mathrm{Aut}(C)$ consists only of two elements, the identity map $\mathrm{id}_C$ and its additive inverse $-\mathrm{id}_C$ (Remark~\ref{rigid:cycle}).
This explains why almost all of our results in Sections~\ref{sec:4}--\ref{sec:6} hold, and stated, for arbitrary (not necessarily dense) subgroups of $\R$.

Basic properties of the subgroup $\inv{G}$ of $\R^\times$ specific to the one-dimensional case $V=\R$ are listed in Proposition~\ref{properties_of_G*}.
In Section~\ref{sec:5} we compute the group $\inv{G}$ for many subgroups $G$ of $\R$. 

Problem~\ref{pro:2} in the one-dimensional case $V=\R$ is studied in Section~\ref{sec:6}. Let us highlight here some connections of our results with number theory.
We prove that, for every transcendental number $x>1$, 
the symmetric cyclic subgroup $A_x = \{\pm x^m: m \in \mathbb Z\}$ of $\R^\times$ coincides with $\inv{\Z[x^{-1},x]}$,
where $\Z[x^{-1},x]$ 
is the ring
of Laurent polynomials with integer coefficients (Proposition~\ref{transcendental:proposition}).
Whether $A_x$ can be represented as $\inv{G}$ for some subgroup $G$ of $\R$ in case $x$ is an algebraic number, is not completely clear (Question~\ref{question_discrete_group}). 
Nevertheless, for an integer number $x>1$, one can find a subgroup $G$ of $\R$ satisfying $A_x=\inv{G}$ if and only if $x$ is a prime number (Proposition~\ref{prop:6.7}).

The finite-dimensional case $V=\R^n$ for $n\ge 2$ is handled in 
Sections~\ref{sec:7}--\ref{sec:9}. 
In this case the group $HGL(\R^n)$ coincides with the general linear group $GL(n,\R)$ (Remark~\ref{linear:maps:in:dimension:n}). 
Another specific result about the finite-dimensional case $V=\R^n$ is 
Theorem~\ref{cardinality:theorem} establishing that $|\mathrm{Aut}(G)|=|\inv{G}|\le |G|$ for every dense subgroup $G$ of $\R^n$ (for all integers $n\ge 1$). In Section~\ref{sec:7} we show that this inequality is no longer valid in the infinite-dimensional case by 
giving an example of a countable dense subgroup $G$ of the complete separable vector space
$\R^{\mathbb{N}}$ whose automorphism group $\mathrm{Aut}(G)$ has cardinality continuum
(Example~\ref{example:in:R^N}).

We also investigate rectangular dense subgroups (that is, dense subgroups of the form $G_1\times G_2\times \dots \times G_n$) of $\R^n$.
For example,  we prove in Proposition~\ref{properties_of_G*_n-b} that a subgroup $H$ of $\R$ is divisible
if and only if $\inv{H^n}$ contains the general linear group $GL(n,\Q)$ 
of all non-degenerated $(n\times n)$-matrices with rational coefficients. On the other hand, in Remark~\ref{8.2} we find two dense divisible subgroups $G_1$, $G_2$ of $\R$ such that 
$\inv{G_1\times G_2}$ does not contain $GL(2,\Q)$.
(Note
that $G_1\times G_2$ is a dense divisible subgroup of $\R^2$.)
In Example~\ref{non-rigid:square} we exhibit a rigid dense subgroup $G$ of $\R$ such that the automorphism group $\mathrm{Aut}(G^2)$ of its square $G^2$ is infinite, so
$G^2$ is not rigid.
In Section~\ref{sec:8} we compute $\inv{G}\cong \mathrm{Aut}(G)$ for many examples of rectangular dense subgroups $G$ of $\mathbb R^2$.

Section~\ref{sec:9} studies the inverse Problem~\ref{pro:2} in the case $V=\R^n$ for $n\ge 2$.
The main result here is Theorem~\ref{theorem:9.4}
stating that if a dense subgroup $G$ of $\R^n$ for $n\ge 2$ satisfies $SL(n,\R)\subseteq \inv{G}$, then $G=\R^n$ and $\inv{G}=GL(n,\R)$. Therefore, if $G$ is a proper dense subgroup of $\R^n$ for $n\ge 2$, then $\inv{G}$ cannot contain the special linear group $SL(n,\R)$ (Corollary~\ref{9:5:new}).
It follows that, for $n\ge 2$, $SL^{\pm}(n)=\{A\in GL(n): \det A=\pm 1\}$ is a symmetric subgroup of $GL(n)$ such that $SL^{\pm}(n)\not=\inv{G}$ for any dense subgroup $G$ of $\R^n$ (Corollary~\ref{SL:corollary}). 

The general linear group $GL(n,\R)$ has a natural Euclidean topology coming from $\R^{n^2}$. Therefore, one can discuss the covering dimension $\dim$ of any subgroup of $GL(n,\R)$. 
We consider the set $D_n = \{\dim \inv{G}: G \mbox{ a dense subgroup of } \mathbb R^n \}$ and show in Theorem~\ref{dimension:theorem} that 
$$
\{mn: m=0,1,\dots,n\}\subseteq D_n\subseteq \{0,1,\dots,n^2\}.
$$
It remains unclear if the equality $D_n=\{0,1,\dots,n^2\}$ holds (Problem~\ref{problem:D_n}).

\section{Proof of Theorem~\ref{isomorphism:theorem}}
\label{sec:2}

The following proposition (whose proof can be extracted from the proof of \cite[Lemma 6.1]{DS}) provides a sufficient condition for the restriction map to be a monomorphism.

\begin{proposition}
\label{prop:2.1:a}
If $X$ is a subset of a topological vector space $V$ such that the smallest subgroup $\grp{X}$ of the additive group $(V,+)$ of $V$ is dense in $V$, then the map $r_X: \inv{X}\to \mathrm{Homeo}(X)$ is a monomorphism between the groups $\inv{X}$ and $\mathrm{Homeo}(X)$.
\end{proposition}
Proof.
Assume that $\varphi,\psi\in HGL(V)$ and $r_X(\varphi)=r_X(\psi)$. Then $\varphi\restriction_X=\psi\restriction_X$ by \eqref{eq:2}.
Since both $\varphi$ and $\psi$ are homomorphisms of $(V,+)$ to itself, we conclude that 
$\varphi\restriction_{\grp{X}}=\psi\restriction_{\grp{X}}$.
Since $\grp{X}$ is dense in $V$ and continuous maps $\varphi,\psi$ coincide on the dense subset $\grp{X}$ of $V$, they must coincide on $V$ as well; that is, $\varphi=\psi$.
$\Box$

\medskip
From Propositions~\ref{basic:proposition} and~\ref{prop:2.1:a}, 
we get the following

\begin{corollary}
\label{monomorphism}
If $G$ is dense subgroup (of the additive group $(V,+)$) of a topological vector space $V$, then $r_G:\inv{G}\to \Aut(G)$
is a monomorphism.
\end{corollary}

We shall need the following folklore fact (cf. \cite[Proposition 7.5(ii)]{HM}).

\begin{fact}\label{continuous_homomorphism}
\label{homomorphisms:are:linear}
Every continuous homomorphism $\varphi:(V_1,+)\to (V_2,+)$ between additive groups $(V_1,+)$ and $(V_2,+)$ of topological vector spaces $V_1$ and $V_2$, is a linear mapping.
\end{fact}

\begin{lemma}
\label{linear:extensions}
Let $G$ be a dense subgroup of the additive group of a complete vector space $V$. Then for every continuous homomorphism 
$f:G\to G$ there exists 
a continuous linear map $\varphi: V\to V$
such that $f=\varphi\restriction_G$. 
\end{lemma}
Proof.
Let $\varphi:\tilde{G}\to \tilde{G}$ be the continuous homomorphism 
extending $f$ over the Raikov completion $\tilde{G}$ of $G$ (\cite[Corollary 3.6.17]{AT}).
Since $G$ is dense in $(V,+)$ and the latter group is Raikov complete by our assumption,
we get 
$\tilde{G}= V$ by \cite[Theorem 3.6.14]{AT}. By Fact~\ref{continuous_homomorphism}, the homomorphism $\varphi$ is a linear mapping.
$\Box$

\bigskip
\noindent
{\bf Proof of Theorem~\ref{isomorphism:theorem}.}
The map $r_G:\inv{G}\to \Aut(G)$ is a monomorphism by 
Corollary~\ref{monomorphism}, so it remains only to check its surjectivity.

Let $f\in \Aut(G)$. Then $f:G\to G$ is a continuous homomorphism, so we can apply Lemma~\ref{linear:extensions} to find a continuous linear map $\varphi: V\to V$
such that $f=\varphi\restriction_G$.
Since $f\in\Aut(G)$, the converse $f^{-1}:G\to G$ of $f$ is also a continuous homomorphism, so we can apply Lemma~\ref{linear:extensions} once again to find 
a continuous linear map $\psi: V\to V$
such that $f^{-1}=\psi\restriction_G$. Now note that 
$$
(\psi\circ \varphi)\restriction_G
=
\psi\circ \varphi\restriction_G
=
\psi\circ f
=
\psi\restriction_G\circ f
=
f^{-1}\circ f=\id_G=\id_V\restriction_G
$$
and
$$
(\varphi\circ \psi)\restriction_G
=
\varphi\circ \psi\restriction_G
=
\varphi\circ f^{-1}
=
\varphi\restriction_G\circ f
=
f\circ f^{-1}=\id_G=\id_V\restriction_G.
$$
Since $G$ is dense in $V$ and both $\varphi$ and $\psi$ are continuous maps, we conclude that
$\psi\circ\varphi=\varphi\circ \psi= \id_V$.
Therefore, $\psi=\varphi^{-1}$.
We have shown that both $\varphi$ and $\varphi^{-1}$ are continuous linear maps from $V$ to $V$, which means that $\varphi\in HGL(V)$.

Note that $\varphi(G)=\varphi\restriction_G(G)=f(G)=G$, as $f\in\Aut(G)$. This means that $\varphi\in\inv{G}$ by \eqref{eq:1}.
Finally,
$r_G(\varphi)=\varphi\restriction_G=f$.
This finishes the proof of surjectivity of $r_G$.
$\Box$

\section{Basic properties of the group $\inv{G}$}
\label{sec:3}

\begin{proposition}
\label{symmetric:proposition}
If $X$ is a subset of a topological vector space $V$ such that
$-X=X$, then $-\inv{X}=\inv{X}$.
\end{proposition}
Proof.
Let $\varphi\in \inv{X}$. 
By \eqref{eq:1}, $\varphi\in HGL(V)$ and $\varphi(X)=X$.
Clearly, $-\varphi\in HGL(V)$ as well. Furthermore,
$-\varphi(X)=-X=X$ by our assumption. This shows that $-\varphi\in \inv{X}$. Since $\varphi\in \inv{X}$ was chosen arbitrarily, this proves the inclusion $-\inv{X}\subseteq \inv{X}$.
From this inclusion we obtain the reverse inclusion 
$\inv{X}=-(-\inv{X})\subseteq -\inv{X}$ as well.
$\Box$

The following corollary provides two necessary conditions on $\inv{G}$:

\begin{corollary}
\label{elementary:properties:of:Phi(G):corollary}
The following two conditions hold for every subgroup $G$ of a topological vector space $V$:
\begin{itemize}
\item[(i)] $-\inv{G}=\inv{G}$;
\item[(ii)] $\{\id_V, -\id_V\}$ is a subgroup of  $\inv{G}$.
\end{itemize}
\end{corollary}
Proof.
Since $-G=G$, item (i) follows from Proposition \ref{symmetric:proposition}. Since $\id_V\in HGL(V)$
and 
$\id_V(G)=G$, we have $\id_V\in \inv{G}$ by \eqref{eq:1}.
Now 
$-\id_V\in \inv{G}$ by item (i).
Finally, $\{\id_V, -\id_V\}$ is clearly a subgroup of $GL(V)$.
$\Box$

\medskip
Recall that an abelian group $G$ is called {\it divisible\/} if for every $x \in G$ and every integer $m \ne 0$ there exists $y \in G$ such that $x = m \cdot y$.
One easily sees that an abelian group $G$ is {\em divisible\/} if and only if $nG=G$ for every positive integer $n$, where $nG=\{ng:g\in G\}$.

A stronger version of item (ii) of Corollary~\ref{elementary:properties:of:Phi(G):corollary} characterizes divisible subgroups of $V$.

\begin{proposition}
\label{divisible:proposition}
  For a subgroup $G$ of a topological vector space $V$, the following conditions are equivalent:
\begin{itemize}
\item[(i)] $G$ is divisible:
\item[(ii)]
$(\mathbb{Q}\setminus\{0\})\cdot \id_V
\subseteq \inv{G}$;
\item[(iii)] 
$(\mathbb{Q}\setminus\{0\})\cdot \inv{G}=\inv{G}$.
\end{itemize}
\end{proposition}
Proof.
(i)$\to$(iii) Since $G$ is a subgroup of a vector space (over $\R$),
one can easily see that the divisibility of $G$ is equivalent to the equality
$(\mathbb{Q}\setminus\{0\})\cdot G=G$.

Suppose now that $\varphi\in\inv{G}$ and $q\in \mathbb{Q}\setminus\{0\}$. 
Then $qG=G$, as $G$ is divisible.
Since $\varphi\in\inv{G}$, $\varphi$ is a linear map, and so
$q\cdot \varphi(G)=\varphi(q\cdot G)=\varphi(G)=G$
by \eqref{eq:1}. From this and \eqref{eq:1}, we conclude that
$q\cdot\varphi\in\inv{G}$. This establishes the inclusion
$(\mathbb{Q}\setminus\{0\})\cdot \inv{G}\subseteq \inv{G}$.
The converse inclusion $\inv{G}\subseteq (\mathbb{Q}\setminus\{0\})\cdot \inv{G}$ is trivial.

(iii)$\to$(ii) follows from $\id_V\in\inv{G}$ given by Corollary~\ref{elementary:properties:of:Phi(G):corollary}(ii).

(ii)$\to$(i)
Suppose that $(\mathbb{Q}\setminus\{0\})\cdot \id_V
\subseteq \inv{G}$. Fix a positive integer $n$ and note that 
$n\cdot \id_V\in \inv{G}$ by our assumption, so
$G=n\cdot \id_V(G)=nG$ by \eqref{eq:1}. This shows that $G$ is divisible.
$\Box$

\begin{proposition}
\label{prop:3.4}
Let $G$ be a 
subgroup of a topological vector space $V$ and
$\psi\in HGL(V)$. Then:
\begin{itemize}
\item[(i)]
$\psi(G)$ is a subgroup of $V$;
\item[(ii)] if $G$ is dense in $V$, then so is $\psi(G)$;
\item[(iii)] if $\varphi\circ\psi=\psi\circ\varphi$ for all $\varphi\in \inv{G}$, then $\inv{G}\subseteq \inv{\psi(G)}$.
\end{itemize}
\end{proposition}
Proof.
From $\psi\in HGL(V)$ and the definition of $HGL(V)$, we get
that $\psi:V\to V$ is a continuous surjective homomorphism (of the additive group $V$ onto itself). 

(i) 
Since $G$ is a subgroup of  $V$ and $\psi$ is a homomorphism, $\psi(G)$ is a subgroup of $V$. 

(ii)
Since $G$ is dense in $V$ and $\psi$ is a continuous surjection,
$\psi(G)$ is dense in $\psi(V)=V$.

(iii) Let $\varphi\in \inv{G}$. Then $\varphi\in HGL(V)$ and $\varphi(G)=G$ by \eqref{eq:1}, so
$$
\varphi(\psi(G))=\varphi\circ\psi(G)=\psi\circ \varphi(G)=\psi(\varphi(G))=\psi(G).
$$
From this, $\varphi\in HGL(V)$ and \eqref{eq:1}, we get $\varphi\in \inv{\psi(G)}$. Since $\varphi\in \inv{G}$ was chosen arbitrarily, this establishes the inclusion $\inv{G}\subseteq \inv{\psi(G)}$.
$\Box$

\section{Automorphism groups of subgroups of $\R$: Basic properties}
\label{sec:4}

\begin{definition}
A topological group $G$ 
satisfying $\Aut(G)=\{\id_G,-\id_G\}$
will be called {\em rigid\/}.
\end{definition}

\begin{remark}
\label{rigid:cycle}
If $G$ is a cyclic subgroup of $\R$, then one easily sees that 
$\Aut(G)=\{\id_G,-\id_G\}$ and $\inv{G}=\{\id_\R,-\id_\R\}$.
This implies that all {\em cyclic subgroups of $\R$ are rigid\/}.
\end{remark}

Additional (non-cyclic) examples of rigid subgroups of $\R$ will be given in Example~\ref{rigid_subgroups_of_R}.

\begin{remark}
\label{non-cyclic:remark}
By \cite[Proposition 19]{M}, {\em every non-cyclic subgroup $G$ of $\R$ is dense in $\R$, so $\Aut(G)$ and $\inv{G}$ are isomorphic\/}
by Theorem~\ref{isomorphism:theorem}. Therefore, we can concentrate on computation of $\inv{G}$.
\end{remark}

\begin{remark}
\label{linear:maps:in:dimension:1}
Every linear map from $\R$ to $\R$ has the form $f_r(x)=r x$ ($x\in \R$) for some $r\in\R$. Since $f_r\circ f_s=f_{rs}$ for $r,s\in \R$, 
the group $GL(1,\R)\cong GL(\R)$ is 
isomorphic to the group $\mathbb R^\times = \mathbb R \setminus \{0\}$ of real numbers with the multiplication as the group operation. Below we will identify these groups.
\end{remark}

For $x\in\R$ and $A\subseteq \R$, we let 
$x\cdot A=A \cdot x=\{a\cdot x:a\in A\}$. 

The subgroup $\{-1, 1\}$ of $\mathbb R^\times$ will be denoted by $\mathbb E$.

\begin{proposition}
\label{properties_of_G*}
\label{property of xG}
Let $G$ be a non-zero subgroup of $\R$. Then the following holds:
\begin{itemize}
\item[(a)] $\inv{G}=\{r\in\R^\times: r\cdot G=G\}$;
\item[(b)]
$\inv{G}\supseteq \mathbb E$;
	\item[(c)]  $\inv{G} \subseteq \left\{\frac{x}{y}: x, y \in G \setminus \{0\}\right\}$;
	\item[(d)] if $1 \in G$, then $\inv{G} \subseteq G \setminus \{0\}$;
\item[(e)] $\inv{x\cdot G} = \inv{G}$ for each $x \in \mathbb R \setminus \{0\}$;	
\item[(f)] $G$ is divisible if and only if $\inv{G} \supseteq \mathbb Q^\times.$
\end{itemize}
\end{proposition}
Proof.
(a) follows from Remark~\ref{linear:maps:in:dimension:1} and~\eqref{eq:1}.

(b) follows from Corollary~\ref{elementary:properties:of:Phi(G):corollary}.

(c) Suppose that $r\in\inv{G}$. Then $r\not=0$ by (a). Since $G$ is non-zero,
there exists $y\in G\setminus\{0\}$. Since $r\in \inv{G}$, from (a) we have $r\cdot y \in G$. Since both $r$ and $y$ are non-zero,
$x=ry\not=0$, so $x\in G\setminus\{0\}$. Now $r=\frac{x}{y}$ with $x,y\in G\setminus\{0\}$.

(d) Let $r\in\inv{G}$. Then $r\not=0$ by (a). Since $1\in G$ by our assumption, it follows from (a) that $r\cdot 1=r\in G$.

(e) Since $\R^\times=GL(\R)=HGL(\R)$ is commutative, applying Proposition~\ref{prop:3.4}(iii) twice, we get 
$$
\inv{G}\subseteq \inv{x\cdot G}\subseteq \invbrackets{\frac{1}{x}\cdot x\cdot G}=\inv{G}.
$$

(f) follows from Proposition~\ref{divisible:proposition}.
$\Box$

\medskip
It follows from Remark~\ref{non-cyclic:remark} and
Proposition~\ref{property of xG}(a) 
that, for a non-cyclic subgroup $G$ of $\R$, 
the problem of computing $\Aut(G)$ is equivalent to a purely algebraic problem of
describing the subgroup 
$\inv{G}$ of $\mathbb R^\times$ consisting of those elements $r \in \mathbb R^\times$ such that $r\cdot G=G$.

\begin{corollary}
\label{cardinality:of:automorphism:group:dimension:1}
  For 
every non-zero
subgroup $G$ of $\mathbb R$, we have either $|\inv{G}| = 2$ or  $\aleph_0 \leq |\inv{G}| \leq |G|$. 
\end{corollary}
Proof. By Proposition~\ref{property of xG}(b), 
$\mathbb E\subseteq \inv{G}$. If $\mathbb E=\inv{G}$, then 
$|\inv{G}| = 2$ holds. Suppose now that $\inv{G}\setminus \mathbb E\not=\emptyset$, and let $r\in \inv{G}\setminus \mathbb E$. Clearly, $r\not=1$. Since $\inv{G}$ is a subgroup of $\R^\times$, we conclude that 
$\{r^n:n=1,2,\dots\}$ is an infinite subset of $\inv{G}$, which shows that $\aleph_0\le |\inv{G}|$.
Finally, $|\inv{G}|\le |(G\setminus\{0\})^2|=|G|$ by Proposition~\ref{property of xG}(c), as $G$ is infinite. 
$\Box$

\section{Automorphism groups of subgroups of $\R$: Examples}
\label{sec:5}

\begin{proposition}\label{properties_of_subringG*} Let $G$ be 
a subring of $\mathbb R$ with $1$.
Then $\inv{G} = U(G) \subseteq G^\times$, where $U(G)$ denotes the subgroup of units (invertible elements) of the ring $G$ and $G^\times = G \setminus \{0\}$.
In particular, 
 $\inv{G} = G^\times$ if and only if $G$ is a  subfield of $\mathbb R$. 
\end{proposition}
Proof.
Since $1 \in G$, we have $\inv{G} \subseteq G^\times$ by Proposition~\ref{properties_of_G*}(d).
Let us show the equality  $\inv{G} = U(G)$. 
It is evident that $r\cdot G = G$ for each $r \in U(G)$. So by Proposition~\ref{property of xG}(a) we have
$\inv{G} \supseteq U(G)$. Consider now any $r \in \inv{G}$. Since $r\cdot G = G$ and $1 \in G$, $r$ is a unit of the ring $G$, i. e. $\inv{G} \subseteq U(G)$.
$\Box$

\begin{corollary}
\label{Q:and:R}
$\inv{\R}=\R^\times$ and $\inv{\mathbb{Q}}=\mathbb{Q}^\times$, where $\mathbb{Q}$ is the additive group of rational numbers and $\mathbb Q^\times = \mathbb Q \setminus \{0\}$ is a subgroup of $\R^\times$. 
\end{corollary}

From Corollary~\ref{Q:and:R} and Proposition~\ref{property of xG}(e), one gets the following

\begin{corollary} Let $G = x\cdot \mathbb Q$, where $x$ is an irrational number. Then  $\inv{G} = \mathbb Q^\times$.  In particular, $G \cap \inv{G} = \emptyset$. $\Box$
\end{corollary}

For 
$A, B \subseteq \mathbb R$, we let 
$A+B=\{a+b: a\in A, b \in B\}$. Clearly, $A+B=B+A$ always holds.

Note that if $G_1$, $G_2$ are subgroups of $\mathbb R$ and $x \in \mathbb R$, then $G_1 \cdot x$ and $G_1 + G_2$ are subgroups of $\mathbb R$. Moreover, if the group $G_1$ is dense in $\mathbb R$, then the group  $G_1 + G_2$ is dense in $\mathbb R$ too.

\begin{example}\label{rigid_subgroups_of_R} 
All subgroups $G$ of $\mathbb R$ in this example (are dense in $\R$ and) satisfy
$\inv{G} = \mathbb E$, so they are rigid.
\begin{itemize}
	\item[(a)] Let $G = \mathbb Z + \mathbb Q \cdot x$, where $x$ is an irrational number such that $x^2 \in \mathbb Q$.
Since $G$ contains $1$, we have $\inv{G}\subseteq G\setminus\{0\}$ by Proposition~\ref{properties_of_G*}(d).

Consider an arbitrary element $z_1 + q_1 \cdot x \in \inv{G}$,
where $z_1\in\Z$ and $q_1\in\Q$. 
Assume that $q_1 \ne 0$. 
Since $x^2\in\Q\setminus\{0\}$,
we can choose $q_2\in\Q$ such that $q_1  \cdot q_2 \cdot x^2 = \frac{1}{2}$. Let $z_2\in\Z$ be arbitrary.
Then $z_2+q_2\cdot x\in G$, and since $z_1 + q_1 \cdot x \in \inv{G}$, we must have
$$
(z_1+q_1 \cdot x) \cdot (z_2+q_2  \cdot x) = (z_1 \cdot z_2 + q_1 \cdot q_2 \cdot x^2) + (z_1 \cdot q_2 + q_1 \cdot z_2) \cdot x\in G;
$$
in particular,
$z_1 \cdot z_2 + q_1 \cdot q_2 \cdot x^2\in\Z$.
On the other hand, $z_1\cdot z_2\in \Z$ and $q_1  \cdot q_2 \cdot x^2 = \frac{1}{2}$ by our choice, which implies $z_1 \cdot z_2 + q_1 \cdot q_2 \cdot x^2\notin\Z$. This contradiction shows that 
$q_1 = 0$, which in turn implies
that $\inv{G} \subseteq \mathbb Z$. Now it is easy to see that $\inv{G} = \mathbb E \subseteq  G$.
	
	\item[(b)] Let $G = \mathbb Q + \mathbb Z \cdot x$, where $x$ is an irrational number such that $x^2 \in \mathbb Q$.
Observe that 
$$\mathbb Z + \frac{1}{x} \cdot \mathbb Q = \mathbb Z + x \cdot \mathbb Q
\mbox{ and }
G = x \cdot \left(\mathbb Z + \frac{1}{x} \cdot \mathbb Q\right)
=
x\cdot (\mathbb Z + x \cdot \mathbb Q),
$$
so
$$
\inv{G}=\inv{x\cdot (\mathbb Z + x \cdot \mathbb Q)}
=
\inv{\mathbb Z + x \cdot \mathbb Q}=\mathbb E 
$$
by
Proposition~\ref{property of xG}(e) and item (a) of this example.
Note that $\inv{G} \subseteq  G$.
	
\item[(c)]	Let $G = \mathbb Z \cdot x+ \mathbb Q \cdot y$, where $x, y$ are irrational numbers such that $x^2, y^2 \in \mathbb Q$ and $x \cdot y \notin \mathbb Q$. 
Observe that 
$$G = x \cdot \left(\mathbb Z + \frac{y}{x} \cdot \mathbb Q\right) = x \cdot (\mathbb Z + \mathbb Q \cdot (x \cdot y)),
$$
so we have $\inv{G} = \mathbb E$ by Proposition~\ref{property of xG}(e) and item (a) of this example (applicable because $(x\cdot y)^2=x^2\cdot y^2\in\Q$ by our assumption on $x$ and $y$). 
Note that $\inv{G} \cap G = \emptyset$.
\end{itemize}
\end{example}
\begin{question}
 Let $G$ be the group $\mathbb Z + \mathbb Q \cdot x$ or $\mathbb Q + \mathbb Z \cdot x$, where $x$ is an irrational number such that $x^2 \notin \mathbb Q$. What is the group $\inv{G}$?
\end{question}	

\begin{example}\label{example_dense_subgring_of_R}  Let $G = \mathbb Q + \mathbb Q \cdot x$, where $x$ is an irrational number such that $x^2 \in \mathbb Q$. 
It is easy to see that $G$ is a subfield of $\mathbb R$. Hence by Proposition~\ref{properties_of_subringG*} we have $\inv{G} = G^\times \supsetneq \mathbb Q^\times$.	
	\end{example}

Our next proposition
generalizes  Example~\ref{example_dense_subgring_of_R}.
\begin{proposition}  Let $G = \mathbb Q + \mathbb Q \cdot x$, where $x$ is an irrational number. 
\begin{itemize}
\item[(i)]
If $x^2 \in G$, then $G$  is the  smallest subfield of $\mathbb R$ containing $x$ and $\inv{G} = G^\times$.
\item[(ii)]
If $x^2\notin G$, then $\inv{G} = \mathbb Q^\times$.
\end{itemize}
\end{proposition}
Proof. Note that $1 \in G$ and $G$ is divisible. So $\mathbb Q^\times \subseteq \inv{G} \subseteq G \setminus \{0\}$
by Proposition~\ref{properties_of_G*} (d), (f). Assume that $\inv{G} \ne \mathbb Q^\times$, i. e. there exists an element $g$ of $\inv{G}$ such that $g = a+ x \cdot b$, where $a, b \in \mathbb Q$ and $b \ne 0$. Since $g \cdot G = G$, we get that $g \cdot x = (a+ x \cdot b) \cdot x \in G$. So $b \cdot x^2 \in G$ and hence $x^2 \in G$.  In particular, $x$ is algebraic and $G = \mathbb Q(x)$  is the  smallest subfield of $\mathbb R$ containing $x$. It follows from Proposition~\ref{properties_of_subringG*} that  $\inv{G} = G \setminus \{0\}=G^\times$. 
$\Box$

\begin{example}\label{dense_subgroups_of_R}  Let $G = \mathbb Q \cdot x+ \mathbb Q \cdot y$, where $x, y$ are irrational numbers such that $x^2, y^2 \in \mathbb Q$ and $x \cdot y \notin \mathbb Q$. Note that $G = x \cdot (\mathbb Q + \frac{y}{x} \cdot \mathbb Q)= x \cdot (\mathbb Q + \mathbb Q \cdot (x \cdot y))$,
so $\inv{G} = (\mathbb Q + \mathbb Q \cdot (x \cdot y))^\times$
by Proposition~\ref{property of xG}(e) and Example~\ref{example_dense_subgring_of_R}. 

We are going to prove that 
$G \cap \inv{G} = \emptyset$.
This equality follows from Claim~\ref{lemma_ortogonality} below.
\begin{claim}
\label{claim:1}
 The equation $p + q\cdot x = r \cdot y$ with the respect to the rational variables $p,q,r$ has only the trivial solution
$p= q = r = 0$.
\end{claim}
Proof. Let $p + q\cdot x = r \cdot y$ be valid for some rational  $p,q,r$. Then $(p + q\cdot x)^2 = (r\cdot y)^2$ holds.
So $p^2 +2 \cdot p \cdot q \cdot x + q^2 \cdot x^2 = r^2 \cdot y^2$, which implies 
$2 \cdot p \cdot q \cdot x = r^2 \cdot y^2 - p^2 - q^2 \cdot x^2 \in \mathbb Q$. If $p \cdot q \ne 0$, then $x\in\Q$, in contradiction with the choice of $x$.
Otherwise, evidently, $p= q = r = 0$.
$\Box$

\begin{claim}\label{lemma_ortogonality} 
The equation 
\begin{equation}
\label{eq:*}
a + b\cdot (x\cdot y) = c \cdot x + d \cdot y
\end{equation}
 with the respect to the rational variables $a,b,c, d$ has
only the trivial solution $a = b = c = d = 0$. 
\end{claim}
Proof. Consider the equation \eqref{eq:*}. If $b = d = 0$, then evidently $a = c = 0.$ Assume now that either $b$ or $d$ is not equal to $0$.
Then 
$$
1 \cdot y = \frac{a - c \cdot x}{d - b \cdot x} = \frac{(a - c \cdot x) \cdot (d + b \cdot x)}{(d - b \cdot x) \cdot (d + b \cdot x)} =
\frac{a\cdot d - c \cdot b \cdot x^2 }{d^2-b^2 \cdot x^2} + \frac{a \cdot b - c \cdot d}{d^2-b^2 \cdot x^2} \cdot x \in \mathbb Q + \mathbb Q \cdot x,
$$
which contradicts Claim~\ref{claim:1}.
$\Box$
\end{example}

\section{Automorphism groups of subgroups of $\R$: Problem~\ref{pro:2}}
\label{sec:6}

\begin{definition}
\label{def:6:1}
For a subgroup $G$ 
of $\mathbb R^\times$, we denote by  $\widehat{G}=\grp{G}$  the subgroup of $\mathbb R$ generated by the set $G$.
\end{definition}

\begin{remark}
\label{remark:6:2}
For a subgroup $G$ 
of $\mathbb R^\times$, 
$$
\widehat{G} = \left\{\sum_{i=1}^n m_i \cdot g_i: m_i \in \mathbb Z, g_i \in G, n = 1, 2, \dots \right\}
$$
 is a  subring of $\mathbb R$.
\end{remark}

\begin{proposition}\label{extension of a multiplicative group} Let $G$ be a 
subgroup of $\mathbb R^\times$.
Then $\inv{\widehat{G}} = U(\widehat{G}) \supseteq G$.
\end{proposition}
Proof.  
Since $\widehat{G}$ is a subring of $\R$ by Remark~\ref{remark:6:2},
it follows from Proposition~\ref{properties_of_subringG*} that 
$\inv{\widehat{G}} = U(\widehat{G})$. Since $G \subseteq U(\widehat{G})$,  we have $\inv{\widehat{G}} \supseteq G$.
$\Box$

\medskip
The following is a reformulation of Problem~\ref{pro:2} in the one-dimensional case.

\begin{problem}
\label{inverse:problem:in:dimension:1}
 Which subgroups $H$ of the group $\mathbb R^\times$ are of the form  $\inv{G}$ for some (dense)
subgroup $G$ of the group  $\mathbb R$?
\end{problem}

A necessary condition is given in Corollary~\ref{elementary:properties:of:Phi(G):corollary}(i): The subgroup $H$ must be {\em symmetric\/}, i.e. satisfy the equation $-H=H$.

Let $x> 1$ be a
real number.
Note that
$$
A_x = \{\pm x^m: m \in \mathbb Z\}
$$
is a symmetric subgroup
of the group $\mathbb R^\times$
which is a discrete subspace of $\mathbb R \setminus \{0\}$.
One can easily see that 
every non-trivial symmetric subgroup $H$ of $\R^\times$ contains
$A_x$ for a suitable real number $x>1$. 

\begin{proposition}
\label{transcendental:proposition}
Let $x>1$ be a transcendental number. Then: 
\begin{itemize}
\item[(i)]
$\inv{\widehat{A_x}}=A_x$;
\item[(ii)]
$\inv{D(\widehat{A_x})}=(\Q\setminus\{0\})\cdot A_x$, where 
$D(\widehat{A_x})$ is the divisible hull of $\widehat{A_x}$.
\end{itemize}
\end{proposition}
Proof. Let $G=\widehat{A_x}$.

(i)
It follows from Remark~\ref{remark:6:2} that
each non-zero element $g\in G$ has the form 
\begin{equation}
\label{repr:g}
g=\sum_{n\in F} z_n x^n,
\end{equation}
 where $F$ is a finite subset of $\Z$
and $z_n\in\Z\setminus\{0\}$ for each $n\in F$. 
Since $x$ is a transcendental number, 
the representation \eqref{repr:g} is unique, 
and therefore
$G$ coincides with the ring
$\Z[x^{-1},x]$ of Laurent polynomials with integer coefficients, so
$\inv{G}=U(\Z[x^{-1},x])$ by Proposition~\ref{extension of a multiplicative group}.
It remains to note that 
$$
U(\Z[x^{-1},x])=\{U(\Z)\cdot x^n:n\in\Z\}=\{\{-1,1\}\cdot x^n:n\in\Z\}=\{\pm x^n:n\in\Z\}=A_x;
$$
see \cite{L}.

(ii) It follows from the argument in the proof of (i) that each non-zero element $g$ of the group $D(G)$ has a unique representation \eqref{repr:g} in which $z_n\in\Q\setminus\{0\}$.
Therefore, $D(G)$ coincides with the ring
$\Q[x^{-1},x]$ of Laurent polynomials with rational coefficients, so
$\inv{D(G)}=U(\Q[x^{-1},x])$ by Proposition~\ref{extension of a multiplicative group}.
It remains to note that 
$$
U(\Q[x^{-1},x])=\{U(\Q)\cdot x^n:n\in\Z\}=\{(\Q\setminus\{0\})\cdot x^n:n\in\Z\}=(\Q\setminus\{0\})\cdot A_x;
$$
see \cite{L}.
$\Box$

\medskip
It follows from Proposition~\ref{transcendental:proposition}(i) that 
each group $A_x$ for a transcendental number $x>1$ has the form 
$\inv{G}$ for some (dense) subgroup $G$ of $\R$. What happens for algebraic numbers $x>1$ remains unclear.

\begin{question}\label{question_discrete_group} 
For which algebraic numbers $x>1$ one can find a (dense) subgroup $G$ of $\R$ such that $\inv{G}=A_x$?
\end{question}

Our next proposition shows that not every algebraic number works.

\begin{proposition}
\label{prop:6.7}
  Let $x$ be an integer $>1$. Then the group $A_x$ coincides with $\inv{G}$ for some (dense) subgroup $G$ of $\mathbb R$ if and only if $x$ is a prime number.
\end{proposition}

Proof. 
 Let $x$ be a prime number. 
 It is easy to see that $G=\Z\cdot A_x$ is a countable dense subgroup of $\mathbb R$ such that $\inv{G} = A_x$, so the ``if'' part  is proved. Let us note that $\inv{G} \subseteq G$ and $|G \setminus \inv{G}| = \aleph_0$.

To prove the ``only if'' part, assume that $A_x=\inv{G}$ for some subgroup $G$ of $\R$. Since $x\in A_x=\inv{G}$, from Proposition~\ref{properties_of_G*}(a) we get $x\cdot G =G$.
Assume that $x$ is not a prime number and fix a proper divisor $d\ge 2$ of $x$. Since $d$ is a divisor of $x$, one can easily deduce from 
$x\cdot G =G$ that the equality $d\cdot G =G$ holds as well.
Then $d\in \inv{G}$ by Proposition~\ref{properties_of_G*}(a).
Since $\inv{G}=A_x$, we get $d\in A_x$. On the other hand,
$2\le d<x$ and the definition of $A_x$ implies that $d\not\in A_x$. 
$\Box$

\medskip
In view of Proposition~\ref
{extension of a multiplicative group},
the following question has some interest.

\begin{question}
\label{hat:question}
  For which dense symmetric subgroups $G$ of $\mathbb R^\times$ does the equality 
$\inv{\widehat{G}} = G$ hold?
\end{question}
 
For every positive number $x$, 
\begin{equation}
\label{def:Q_x}
\Q_x=(\Q\setminus\{0\})\cdot A_x
\end{equation}
is a dense symmetric subgroup of $\mathbb R^\times$. 

\begin{proposition}\label{dense_transcendental} If $x$ is a transcendental number,
then $\inv{\widehat{\Q_x}} = \Q_x$. 
	
\end{proposition}
Proof. Note that $\widehat{\Q_x}=D(\widehat{A_x})$, where $D(\widehat{A_x})$ is the divisible hull of $\widehat{A_x}$.
Now the conclusion follows from Proposition~\ref{transcendental:proposition}(ii) and~\eqref{def:Q_x}.
 $\Box$

\medskip
The next example demonstrates that the above proposition does not hold for algebraic numbers $x$, thereby establishing that the equality in Question~\ref{hat:question} can fail for some dense symmetric subgroups $G$ of $\R^\times$.

\begin{example}
\label{example:of:non-coincidence}
 Let $x$ be an irrational number such that $x^2 \in \mathbb Q$. 
Then $\Q_x=(\Q\setminus\{0\})\cdot \{1,x\}$, and so
$\widehat{\Q_x} = \mathbb Q + \mathbb Q \cdot x$ by Definition~\ref{def:6:1}.
From Example~\ref{example_dense_subgring_of_R}, we get 
	  $\inv{\widehat{\Q_x}} = (\mathbb Q + \mathbb Q \cdot x)^\times$. Note that 
$(\mathbb Q + \mathbb Q \cdot x)^\times\supsetneq \Q_x$, as
$1+x \in (\mathbb Q + \mathbb Q \cdot x)^\times \setminus \Q_x$.	  
\end{example}

\begin{proposition}\label{dense_subfield} If $G= F^\times$, where $F$ is a (dense) subfield of $\mathbb R$, then $G$ is a (dense) symmetric subgroup of $\mathbb R^\times$ satisfying $\inv{\widehat{G}} = G$.
\end{proposition}
Proof. It is clear that $F^\times$ is a dense symmetric subgroup of $\mathbb R^\times$ and  $\widehat{G} = F$. Hence by  Proposition~\ref{properties_of_subringG*} we get 
$\inv{\widehat{G}} = \widehat{G}^\times = F^\times = G$.
$\Box$

\begin{remark}
For a symmetric subgroup $G=\{-1,1\}$ of $\R^\times$, one has 
$\widehat{G}=\Z$, and so $\inv{\widehat{G}}=\inv{\Z}=\{-1,1\}=G$
by Remark~\ref{rigid:cycle}.
\end{remark}

\section{Automorphism groups of dense subgroups of $\mathbb R^n$ for $n\ge 2$: Basic properties}
\label{sec:7}

\begin{lemma}
\label{dense:subspace:basis:lemma}
Every dense subset of $\R^n$ (for $n\ge 1$) contains $n$ vectors forming a basis of the vector space $\R^n$.
\end{lemma}

The following proposition extends the second inequality in Corollary~\ref{cardinality:of:automorphism:group:dimension:1} to the finite-dimensional case.

\begin{theorem}
\label{cardinality:theorem}
If $G$ is a dense subgroup of $\R^n$ for some integer $n\ge 1$, then
$|\inv{G}|\le |G|$ and $|\mathrm{Aut}(G)|\le |G|$. 
\end{theorem}
Proof.
Use Lemma~\ref{dense:subspace:basis:lemma} to
find $n$ vectors
$\overline{x}_1, \overline{x}_2,\dots,\overline{x}_n\in G$ which form a basis of the vector space $\R^n$.

Let $\varphi\in \inv{G}\subseteq GL(\R^n)$ be arbitrary. Then $\varphi$ is uniquely determined by $n$ vectors 
$\varphi(\overline{x}_1)$, $\varphi(\overline{x}_2),\dots,\varphi(\overline{x}_n)
$
each of which belongs to $G$ by \eqref{eq:1}. Therefore, the map 
which sends $\varphi\in \inv{G}$ to the element $(\varphi(\overline{x}_1), \varphi(\overline{x}_2),\dots,\varphi(\overline{x}_n))$ of $G^n$ is injective. 
This shows that $|\inv{G}|\le |G^n|$. Since $G$ is dense in $\R^n$, 
it must be infinite, which implies $|G^n|=|G|$. Finally, $|\mathrm{Aut}(G)|=|\inv{G}|$ by Theorem~\ref{isomorphism:theorem}. 
$\Box$

\medskip
The following example shows that this theorem no longer holds in the infinite-dimensional case, even for dense subgroups of complete separable metric vector spaces.

\begin{example}
\label{example:in:R^N}
Consider the 
vector space $\R^{\mathbb{N}}$ equipped with the Tychonoff product topology and its countable dense subgroup
$$
G=\{(x_n)_{n\in\mathbb{N}}: x_n\in\Q\text{ for all }n\in\mathbb{N},
\text{ and }
x_n=0
\text{ for all but finitely many }
n\}.
$$
Then 
$|\mathrm{Aut}(G)|=|\inv{G}|\ge 2^{\aleph_0}>\aleph_0=|G|$.
Indeed,
let $S(\mathbb{N})$ denote the group of all permutations 
of the set $\mathbb{N}$ of natural numbers. Every permutation $f\in S(\mathbb{N})$ generates an element $\varphi_f\in HGL(\R^{\mathbb{N}})$, defined by $\varphi_f((x_n)_{n\in\mathbb{N}})=
(x_{f(n)})_{n\in\mathbb{N}}$ for each $(x_n)_{n\in\mathbb{N}}\in\R^{\mathbb{N}}$, such that 
$\varphi_f(G)=G$. Therefore, $\{\varphi_f: f\in S(\mathbb{N})\}\subseteq \inv{G}$ by \eqref{eq:1}.
Since the map $f\mapsto\varphi_f$ is an injection of 
$S(\mathbb{N})$ into $HGL(\R^{\mathbb{N}})$, 
this implies
$2^{\aleph_0}=|S(\mathbb{N})|=|\{\varphi_f: f\in S(\mathbb{N})\}|\le |\inv{G}|$. Finally, recall that $|\mathrm{Aut}(G)|=|\inv{G}|$ by Theorem~\ref{isomorphism:theorem}.
\end{example}

\begin{remark}
\label{linear:maps:in:dimension:n}
Let $n\ge 2$ be an integer.
Every linear map from $\R^n$ to $\R^n$ has the form $\varphi_A(\overline{x})=\overline{x}\cdot A$ ($\overline{x}\in \R^n$) for some $(n\times n)$-matrix $A$. Since $\varphi_B\circ \varphi_A=\varphi_{AB}$ for $A,B\in GL(n,\R)$, 
the group $GL(\R^n)$ is 
isomorphic to the group $GL(n,\R)$ of $(n\times n)$-matrices with the multiplication of matrices as the group operation. Below we will identify these groups. 
\end{remark}

For simplicity, we shall write
$GL(n)$ instead of $GL(n,\R)$.

For a subset $X$ of $\R^n$ and a matrix $A\in GL(n)$, we 
let $X\cdot A=\{\overline{x}\cdot A: \overline{x}\in X\}$, where 
$\overline{x}\cdot A$ is the matrix multiplication of the row vector $\overline{x}$ and the matrix $A$.

The following result is an analogue of items (a) and (e) of Proposition~\ref{properties_of_G*}.

\begin{proposition}\label{property of xG n > 1} Let $G$ be a 
subgroup of $\mathbb R^n$ for an integer $n > 1$.
Then:
\begin{itemize}
	\item[(a)] $\inv{G}=\{A \in GL(n): G \cdot A = G\}$.
	\item[(b)]  
For every $x \in \mathbb R \setminus \{0\}$, we have  
	$\inv{G \cdot xI} = \inv{G}$, where
$I$ is the identity matrix of $GL(n)$. 
\end{itemize}
\end{proposition}
Proof. (a) follows from Remark~\ref{linear:maps:in:dimension:n} and~\eqref{eq:1}.

(b)
Note that the scalar matrix $\lambda I$ commutes with every matrix in $GL(n)$, so from Remark~\ref{linear:maps:in:dimension:n} and Proposition~\ref{prop:3.4}(iii), we get 
$$
\inv{G}\subseteq \inv{G\cdot xI}\subseteq \invbrackets{G\cdot xI\cdot \frac{1}{x}I}=\inv{G\cdot I}=\inv{G}.
\hskip15pt\Box$$

\begin{corollary}
\label{atomorphisms:of:R^n}
$\mathrm{Aut}(\R^n)\cong \inv{\R^n}=GL(n)$.
\end{corollary}
Proof.
The equality $\inv{\R^n}=GL(n)$ follows from item (a) of Proposition~\ref{property of xG n > 1} and the definition of $GL(n)$, while the isomorphism $\mathrm{Aut}(\R^n)\cong \inv{\R^n}$ was proved in Theorem~\ref{isomorphism:theorem}.
$\Box$

\medskip
Note that
\begin{equation}
\label{GL_Q(n)}
GL_\mathbb Q(n) = \{(a_{i,j})_{i, j \leq n}   \in GL(n): a_{i,j} \in \mathbb Q\text{ for all }i,j\le n\}
\end{equation}
is a subgroup of $GL(n)$. In particular, $A\in GL_\mathbb Q(n)$ implies $A^{-1}\in GL_\mathbb Q(n)$.

The following result is a certain multi-dimensional analogue of item (f) of Proposition~\ref{properties_of_G*}.

\begin{proposition}
\label{properties_of_G*_n-b}
Let $n\ge 1$ be an integer.
A subgroup $H$ of $\R$ is divisible if and only if $\inv{H^n} \supseteq GL_\mathbb Q(n)$.
\end{proposition}
Proof.
Suppose that $H$ is divisible.
\begin{claim}
\label{claim:l}
$H^n \cdot B \subseteq H^n$ for every $B\in GL_{\Q}(n)$.
\end{claim}
Proof.
Let $B =  (b_{i,j})_{i, j \leq n} \in GL_\mathbb Q(n)$.
Fix $i, j \leq n$. Then $b_{i,j}\in\Q$ by \eqref{GL_Q(n)}.
If $b_{i,j}=0$, then $H \cdot  b_{i,j} = \{0\}\subseteq H$.
If $b_{i,j} \ne 0$, then $b_{i,j}\in\Q^\times$, so
$b_{i,j} \in \inv{H}$ by Proposition~\ref{properties_of_G*}(f), 
which implies 
$H \cdot b_{i,j} = H$ by Proposition~\ref{property of xG}(a).

Let $\overline{x} = (x_1, \dots, x_n)  \in H^n$. 
Then $x_i\in H$ for every $i\le n$, so from what was proved in the previous paragraph we conclude that 
$\sum_{i=1}^n x_i \cdot b_{i, j}\in H$ for every $j\le n$, which implies
$$
\overline{x} \cdot B = \left(\sum_{i=1}^n x_i \cdot b_{i, 1}, \dots, \sum_{i=1}^n x_i \cdot b_{i, n}\right) \in H^n.
$$
We have proved that $H^n \cdot B \subseteq H^n$.
$\Box$ 

\medskip  
Let $B\in GL_{\Q}(n)$ be arbitrary.
Since 
$B^{-1} \in  GL_\mathbb Q(n)$,
from Claim~\ref{claim:l} we get
$H^n \cdot B^{-1} \subseteq H^n$,
so 
$H_n=H_n\cdot E= (H^n \cdot B^{-1}) \cdot B \subseteq H^n \cdot B$. Combining this with Claim~\ref{claim:l}, we get
$H^n \cdot B =  H^n$.
Therefore, $B \in \inv{H^n}$
by Proposition~\ref{property of xG n > 1}(a).
This establishes the inclusion $B\subseteq \inv{H^n}$.

Suppose now that $H$ is not divisible. Then there exists $0 \ne m \in \mathbb Z$ such that $m\cdot H \ne H$. Let $A = m \cdot I$, where $I$
is the identity $(n \times n)$-matrix. Clearly,
$A \in  GL_\mathbb Q(n)$ and 
$H^n \cdot A \ne H^n$.  
From the latter and Proposition~\ref{property of xG n > 1}(a), 
we get $A\not\in\inv{H^n}$.
$\Box$

\medskip
This proposition no longer holds even for the product of two (distinct) divisible subgroups of $\R$; see Remark~\ref{8.2} below. 

The following result is a certain multi-dimensional analogue of item (c) of Proposition~\ref{properties_of_G*}.

\begin{proposition}\label{properties_of_G*_n} Let $G = G_1 \times \dots \times G_n$, where each $G_i$ is a 
subgroup of $\mathbb R$.
Suppose that $A = (a_{i,j})_{i, j \leq n} \in \inv{G}$,
$i, j \leq n$ and $x \in G_i \setminus \{0\}$.
Then $a_{i, j} = \frac{y}{x}$ for some $ y \in G_j$.
\end{proposition}
Proof. 
Let $i \leq n$  and $x \in G_i \setminus \{0\}$. Note that the vector $\overline{x}_i = (x_1, \dots, x_n)$ such that $x_i = x$ and $x_k = 0$ for all $k \ne i$, belongs to $G$. Since $A \in \inv{G}$, we have  
\begin{equation}
\label{eq:***}
\overline{x}_i  \cdot A = (x \cdot a_{i, 1}, \dots,  x \cdot a_{i, n}) \in G
\end{equation}
 by Proposition~\ref{property of xG n > 1}(a). 
Let $j\le n$.
It follows from \eqref{eq:***} and our assumption on $G$ that 
$x\cdot a_{i,j}\in G_j$.
Thus, $a_{i,j} = \frac{y}{x}$ for some $y \in G_j$.
$\Box$

\begin{corollary}
\label{Aut(Q^n)}
$\mathrm{Aut}(\Q^n)\cong \inv{\Q^n}=GL_{\Q}(n)$ for every integer $n\ge 1$.
\end{corollary}
Proof.
The inclusion $\inv{\Q^n}\subseteq GL_{\Q}(n)$ follows from Proposition~\ref{properties_of_G*_n}. The converse inclusion 
$GL_{\Q}(n)\subseteq \inv{G}$ follows from Proposition~\ref{properties_of_G*_n-b}, as $\Q$ is divisible.
Finally, the isomorphism $\mathrm{Aut}(\Q^n)\cong \inv{\Q^n}$ was proved in Theorem~\ref{isomorphism:theorem}, as $\Q^n$ is dense in $\R^n$.
$\Box$

\section{Automorphism groups of dense subgroups of $\mathbb R^2$: Examples}
\label{sec:8}

Our first example in this section is a particular case of a more general Proposition~\ref{Q^n:times:R^m} which shall be proved in the next section.

\begin{example} 
\label{Q:times:R}
The subgroup $G = \mathbb Q \times \mathbb R$ of $\R^2$
satisfies 
$$
\inv{G} = \left\{ \left(\begin{array}{rr}  a & b \\ 0 & c  \end{array}\right) \in GL(2): a \in \mathbb Q, b, c \in \mathbb R \right\}.
$$
\end{example} 

\begin{example}
\label{triple:example}
All subgroups $G$ of $\mathbb R^2$ in this example are 
dense in $\R^2$.
\begin{itemize}

	\item[(a)] Let $G = (\mathbb Q \cdot x) \times (\mathbb Q \cdot x) = (\mathbb Q \cdot x)^2$, where $x$ is an irrational number.
	Then $\inv{G} = GL_\mathbb Q(2).$ Indeed, $G = x \cdot (\mathbb Q \times \mathbb Q) = (\Q^2) \cdot (xI)$, where $I$ is the identity $(2 \times 2)$-matrix.  
By Proposition~\ref{property of xG n > 1}(b) and Corollary~\ref{Aut(Q^n)}, we have $\inv{G} = \inv{\mathbb Q^2 \cdot xI} = \inv{\mathbb Q^2} = GL_\mathbb Q(2).$
	
	\item[(b)] Let $G = \mathbb Q \times (\mathbb Q \cdot x)$, where $x$ is an irrational number such that $x^2 \in \mathbb Q$. 
Then 
$$
\inv{G} = \left\{ \left(\begin{array}{rr}  a & b\cdot x \\ c \cdot x & d  \end{array}\right) \in GL(2): a, b, c, d \in \mathbb Q \right\}.
$$ 
Indeed, first let us consider some matrix $\left(\begin{array}{rr}  a & b\cdot x \\ c \cdot x & d  \end{array}\right),$ where $a, b, c, d \in \mathbb Q,$  and 
note that 
\begin{itemize}
	\item[(i)] $(q_1, q_2 \cdot x) \cdot \left(\begin{array}{rr}  a &b\cdot x \\ c \cdot x & d  \end{array}\right) = (q_1 \cdot a +q_2 \cdot c \cdot x^2, (q_1 \cdot b + q_2 \cdot d) \cdot x) \in G$ for every $(q_1, q_2 \cdot x)  \in G$, and
	
	\item[(ii)] for any $p, r \in \mathbb Q$, the system 
$$
	\left\{\begin{array}{l}  a \cdot u + (c \cdot x) \cdot v = p  \\  (b \cdot x)\cdot  u + d \cdot v = r \cdot x  \end{array} \right. 
$$
  has  only one solution such that $u = q_1\in \mathbb Q$ and $v = q_2  \cdot x$, where $q_1, q_2 \in \mathbb Q$. 
\end{itemize}

The items (i) and (ii) imply the inclusion
$$
\inv{G} \supseteq \left\{ \left(\begin{array}{rr}  a & b\cdot x \\ c \cdot x & d  \end{array}\right) \in GL(2): a, b, c, d \in \mathbb Q \right\}.
$$
 
 The reverse inclusion   
  
$$
\inv{G} \subseteq \left\{ \left(\begin{array}{rr}  a & b\cdot x \\ c \cdot x & d  \end{array}\right) \in GL(2): a, b, c, d \in \mathbb Q \right\}
$$ 
follows from Proposition~\ref{properties_of_G*_n}.

	\item[(c)] Let $G = (\mathbb Q \cdot x) \times (\mathbb Q \cdot y)$, where $x, y$ are  irrational numbers such that $x^2, y^2 \in \mathbb Q$ and $x \cdot y \notin \mathbb Q$.
	Then 
$$
\inv{G} = \left\{ \left(\begin{array}{rr}  a & b\cdot (x \cdot y) \\ c \cdot (x \cdot y) & d  \end{array}\right) \in GL(2): a, b, c, d \in \mathbb Q \right\}.
$$ 
	
	Indeed, since $x^2\in\Q$ by our assumption,
$$
G = x \cdot \left(\mathbb Q \times \left(\frac{y}{x} \cdot \mathbb Q\right)\right) = 
x\cdot\left(\Q\times((x\cdot y)\cdot\Q)\right)
=
(\Q\times((x\cdot y)\cdot\Q))
\cdot (x \cdot I),
$$
 where $I$ is the identity $(2 \times 2)$-matrix. 
Since $(x\cdot y)^2=x^2\cdot y^2\in \Q$ by our assumption,
from Proposition~\ref{property of xG n > 1}(b) and item (b) of this example,
 we get 
$$
\inv{G}=\inv{(\mathbb Q \times (x \cdot y) \cdot \mathbb Q)\cdot xI}
=\inv{\mathbb Q \times (x \cdot y) \cdot \mathbb Q}
=
$$
$$
  = \left\{ \left(\begin{array}{rr}  a & b\cdot (x \cdot y) \\ c \cdot (x \cdot y) & d  \end{array}\right) \in GL(2): a, b, c, d \in \mathbb Q \right\}.
$$ 
\end{itemize}
\end{example}

\begin{remark}
\label{8.2}
Both $\mathbb Q$ and $\mathbb Q \cdot x$ are divisible,
yet
$\inv{(\Q\times(\mathbb Q \cdot x)}$ does not contain $GL_\mathbb Q(2)$ by item (b) of Example~\ref{triple:example}.
This shows the limit of Proposition~\ref{properties_of_G*_n-b}.	
\end{remark}

Let  
$$
GL_\mathbb Z(n) = \{A=(a_{i,j})_{i, j \leq n} \in GL(n): a_{i,j}\in\Z
\text{ for all }i,j\le n\}
$$
Note that 
$$
E_\mathbb Z(n) = \{A \in GL_\mathbb Z(n): A^{-1} \in GL_\mathbb Z(n)\}
$$
is a subgroup of $GL(n)$.
In particular,
if $A \in E_\mathbb Z(n)$ then $A^{-1} \in E_\mathbb Z(n)$.

\begin{remark} $E_\mathbb Z(1) = \mathbb E$ and $E_\mathbb Z(n) = \{A \in GL_\mathbb Z(n): \det A = \pm 1\}$. The second equality can be seen as follows.
If $A \in E_\mathbb Z(n)$, then $A, A^{-1} \in GL_\mathbb Z(n)$ and hence $\det A, \det A^{-1} \in  \mathbb Z$. Since $\det A \cdot \det A^{-1} = 1$, one has $\det A = \pm 1$, i. e. $E_\mathbb Z(n) \subseteq \{A \in GL_\mathbb Z(n): \det A = \pm 1\}$. If 
$A \in GL_\mathbb Z(n)$ and $\det A = \pm 1$, then by the known formulas for the elements of $A^{-1}$ we have $A^{-1} \in GL_\mathbb Z(n)$. So  $E_\mathbb Z(n) \supseteq \{A \in GL_\mathbb Z(n): \det A = \pm 1\}$. 
\end{remark}

\begin{lemma}\label{the last lemma} $\mathbb Z^n \cdot C = \mathbb Z^n$ for each $C \in E_\mathbb Z(n)$.
\end{lemma}
Proof. It is clear that  $\mathbb Z^n \cdot C \subseteq \mathbb Z^n$. In particular, $(\mathbb Z^n \cdot C) \cdot C^{-1}
= \mathbb Z^n   \subseteq \mathbb Z^n \cdot C^{-1}$. 
$\Box$

\medskip
For sets $A,B\subseteq \R^2$, we let $A+B=\{\overline{a}+\overline{b}: \overline{a}\in A, \overline{b}\in B\}$.

We have seen in Example~\ref{rigid_subgroups_of_R} that the subgroup $A=\mathbb Z + \mathbb Q \cdot x$ of $\R$, for an irrational number $x$ such that  $x^2 \in \mathbb Q$, is rigid.
Our next example shows that the square $A^2$ of this group (is dense in $\R^2$ and) is not rigid, and in fact, has an infinite automorphism group.

\begin{example} 
\label{non-rigid:square}
Let $G = (\mathbb Z + \mathbb Q \cdot x) \times (\mathbb Z + \mathbb Q \cdot x) = (\mathbb Z + \mathbb Q \cdot x)^2,$ where $x$ is an irrational number such that  $x^2 \in \mathbb Q$. Note that $G$ is a dense  subgroup of $\mathbb R^2$ and the dense subgroup $\mathbb Z + \mathbb Q \cdot x$ of $\mathbb R$ is not divisible. Moreover,  $G = \mathbb Z^2 + (\mathbb Q \cdot x)^2$.

Let is check that 
$\inv{G} \supseteq E_{\mathbb Z}(2)$.
Let $C \in E_{\mathbb Z}(2)$ be arbitrary. Then 
$$
G \cdot C = (\mathbb Z + \mathbb Q \cdot x)^2 \cdot C = (\mathbb Z^2 + (\mathbb Q \cdot x)^2) \cdot C = \mathbb Z^2 \cdot C + (\mathbb Q \cdot x)^2 \cdot C. 
$$
Now $\mathbb Z^2 \cdot C=\Z^2$
by Lemma~\ref{the last lemma}. Clearly,
$(\mathbb Q \cdot x)^2 \cdot C=(\mathbb Q \cdot x)^2$.
Therefore,
$G \cdot C = \mathbb Z^2 + (\mathbb Q \cdot x)^2 = G$,
so $C\in\inv{G}$ by Proposition~\ref{property of xG n > 1}(a).
\end{example}

\section{Automorphism groups of dense subgroups of $\mathbb R^n$ for $n\ge 2$: Problem~\ref{pro:2}}
\label{sec:9}

In view of Remark~\ref{linear:maps:in:dimension:n},
the following is a reformulation of Problem~\ref{pro:2} in the finite-dimensional case. (The symmetricity condition is necessary by Corollary~\ref{elementary:properties:of:Phi(G):corollary}(i).)

\begin{problem}
\label{problem:1.4:R^n}
 Let $n\ge 2$ be an integer. 
 Which symmetric subgroups of the group $GL(n)$ can be realized as $\inv{G}$ for some dense
subgroup $G$ of the group  $\mathbb R^n$?
\end{problem}

Since $\R^n$ has many countable dense subgroups $G$ and 
the subgroup $\inv{G}$ for such a $G$ is countable by Theorem~\ref{cardinality:theorem}, 
the countable version of Problem~\ref{problem:1.4:R^n} appears to be interesting as well.

\begin{problem}
 Let $n\ge 2$ be an integer. 
 Which countable symmetric subgroups of the group $GL(n)$ can be realized as $\inv{G}$ for some countable dense
subgroup $G$ of the group $\mathbb R^n$?
\end{problem}

In this section we obtain some partial solutions to Problem~\ref{problem:1.4:R^n}.

\begin{lemma}
\label{semigroup:lemma}
For every positive real number $r$, the smallest subsemigroup $G$ of the additive semigroup $(\R^2,+)$ containing the circle $C_r=\{\overline{x}\in \R^2: |\overline{x}|=r\}$ coincides with $\R^2$.
\end{lemma}
Proof.
Note that $\{(x,0):-2r\le x\le 2r\}\cup\{(0,y):-2r\le y\le 2r\}\subseteq C_r+C_r\subseteq G$. Since $G$ is closed under taking arbitrary finite sums, 
from this one easily concludes that
$\{(x,0):x\in\R\}\cup \{(0,y):y\in\R\}\subseteq G$, which in turn implies
$\R^2\subseteq G$. The converse inclusion $G\subseteq \R^2$ is trivial.
$\Box$

Recall that the subgroup 
$$
SO(2)=\{A\in GL(2): A^t A = A A^t=I: \det A=1\}
$$
is called the {\em special orthogonal group of dimension $2$\/}.
(Here $A^t$ denotes the transpose of  $A$.)

\begin{lemma}
\label{two-dimensional:lemma}
If 
$G$ be a non-trivial subgroup of $\R^2$ satisfying 
$SO(2)\subseteq \inv{G}$,
then 
$G=\R^2$, and so $\inv{G}=GL(2)$.
\end{lemma}
Proof.
Fix a non-zero element $\overline{g_0}$ in $G$. Then $r=|\overline{g_0}|>0$.
Recall that 
$SO(2)$
is the group of rotations of the plane $\R^2$, so from $\overline{g_0}\in G$, 
$SO(2)\subseteq \inv{G}$
and 
Proposition~\ref{property of xG n > 1}(a), we obtain that 
$G$ must contain the circle $C_r$. Hence
$G=\R^2$ by Lemma~\ref{semigroup:lemma}, and so
$\inv{G}=\inv{\R^2}=GL(2)$ by 
Corollary~\ref{atomorphisms:of:R^n}.
$\Box$

\begin{corollary}
The subgroup $SO(2)$ of $GL(2)$ (is symmetric but) cannot be realized as $\inv{G}$ for any non-trivial subgroup $G$ of $\R^2$.
\end{corollary}

Recall that the subgroup
$$
SL(n)=\{A\in GL(n): \mathrm{det} A=1\}
$$
of $GL(n)$ is called the {\em special linear group
of dimension $n$\/}.

\begin{theorem}
\label{theorem:9.4}
Let $n\ge 2$ be an integer and let
$G$ be a dense subgroup of $\R^n$ satisfying 
$SL(n)\subseteq \inv{G}$.
Then $G=\R^n$, and so $\inv{G}=GL(n)$.
\end{theorem}
Proof.
Use Lemma~\ref{dense:subspace:basis:lemma} to
find $n$ vectors
$\overline{x}_1, \overline{x}_2,\dots,\overline{x}_n\in G$ which form a basis of the vector space $\R^n$.
Let $\overline{e}_1, \overline{e}_2,\dots,\overline{e}_n$ be the standard basis of $\R^n$; that is, the $i$th coordinate of the vector $\overline{e}_i$ is equal to $1$ and all its other coordinates are equal to $0$.
Let $\pi$ be the linear map which sends each vector $\overline{x}_i$ to the basic vector $\overline{e}_i$, and let $P\in GL(n)$ be the matrix corresponding 
to this map; that is, $\pi=\varphi_P$. 
Note that $\pi\in GL(n)=GL(\R^n)=HGL(\R^n)$, so
$H=\pi(G)$ is a dense subgroup of $\R^n$ 
by Proposition~\ref{prop:3.4}(i),(ii).
Clearly, 
$\overline{e}_1, \overline{e}_2,\dots,\overline{e}_n\in H$. One can easily see that 
$\inv{G}=P\cdot \inv{H}\cdot P^{-1}$, and so
$\inv{H}=P^{-1}\cdot \inv{G}\cdot P$. 
Note that $A\in SL(n)$ if and only if 
$P^{-1}\cdot A\cdot P\in SL(n)$, as $\det A =\det (P^{-1})\cdot \det A\cdot \det P=
\det(P^{-1} \cdot A\cdot P)$.
Since $SL(n)\subseteq \inv{G}$ by our assumption,
\begin{equation}
\label{equation:SL}
SL(n)=P^{-1}\cdot SL(n)\cdot P\subseteq 
P^{-1}\cdot \inv{G}\cdot P=\inv{H}.
\end{equation}

Fix $i=1,2,\dots,n-1$. Consider the family $\Psi_i$ of all linear maps from $\R^n$ to $\R^n$ which are rotations of the plane passing through the $i$th and $(i+1)$th axis and are the identity elsewhere. Clearly, $\Psi_i\subseteq SL(n)$. Combining this with \eqref{equation:SL}, we get $\Psi_i\subseteq \inv{H}$. Since $\overline{e}_i\in H$, arguing as in the proof of Lemma~\ref{two-dimensional:lemma}, we conclude that the circle
$\{\psi(\overline{e}_i): \psi\in \Psi_i\}$ lies inside of $H$, and from this and Lemma~\ref{semigroup:lemma} we get
\begin{equation}
\label{eq:10}
\{0\}\times\dots\times\{0\}\times \R\times\R\times \{0\}\times\dots\times\{0\}\subseteq H,
\end{equation}
where two copies of $\R$ are on the $i$th and $(i+1)$th place, respectively.

Since \eqref{eq:10} holds for every $i=1,2,\dots,n-1$ and $H$ is a subgroup of $\R^n$, this gives $\R^n\subseteq H$.
Since the reverse inclusion $H\subseteq \R^n$ trivially holds,
we obtain the equality $H=\R^n$. Since $H=\pi(G)$ and $\pi$ is a bijection,
we get $G=\pi^{-1}(H)=\pi^{-1}(\R^n)=\R^n$.
Therefore, $\inv{G}=\inv{\R^n}=GL(n)$ by 
Corollary~\ref{atomorphisms:of:R^n}.
$\Box$

\begin{corollary}
\label{9:5:new}
Let $G$ be a proper dense subgroup of $\R^n$ for $n\ge 2$.
Then $\inv{G}$ does not contain $SL(n)$.
\end{corollary}

\begin{corollary}
\label{cor:9.5}
Let $H$ be a subgroup of $GL(n)$ satisfying $SL(n)\subseteq H\subsetneq GL(n)$.
Then $H\not=\inv{G}$ 
for any dense subgroup $G$ of $\R^n$. 
\end{corollary}

\begin{corollary}
\label{SL:corollary}
\begin{itemize}
\item[(i)]
If $n\ge 2$ is an even number, then $SL(n)$ is a symmetric subgroup of $GL(n)$ such that $SL(n)\not=\inv{G}$ for any dense subgroup $G$ of $\R^n$. 
\item[(ii)]
$SL^{\pm}(n)=\{A\in GL(n): \det A=\pm 1\}$ is a symmetric subgroup of $GL(n)$ such that $SL^{\pm}(n)\not=\inv{G}$ for any dense subgroup $G$ of $\R^n$. 
\end{itemize}
\end{corollary}

Let $M_{p,q}$ be the set of $(p \times q)$-matrices with real coefficients, and $O_{p,q}$ be the null matrix of $M_{p,q}$.

\begin{proposition}\label{Q^n:times:R^m} 
The subgroup $G = \mathbb Q^n \times \mathbb R^m$ of $\R^{n+m}$, where $n, m  \geq 1$,
satisfies 
$$
\inv{G} = \left\{ \left(\begin{array}{rr}  A & B \\ O & C  \end{array}\right) \in
M_{n+m, n+m}: 
A \in GL_\mathbb Q (n), B \in M_{n, m}, C \in GL (m),  O = O_{m, n} \right\}. $$
\end{proposition}
Proof. We prove the equality in two steps.

First, let us observe that any matrix 
\begin{equation}
\label{triangulation}
D = \left(\begin{array}{rr}  A & B \\ O & C  \end{array}\right) \in M_{n+m, n+m}
\end{equation}
such that 
\begin{equation}
\label{eq:ABCO}
A \in GL_\mathbb Q (n), B \in M_{n, m}, C \in GL (m),  O = O_{m, n}
\end{equation}
  belongs to $\inv{G}$. (Here $O_{m, n}$ denotes the zero $(m\times n)$-matrix.)
Indeed, note that $\det D = \det A \cdot \det C \ne 0$, so $D \in GL(n+m)$. A straightforward check shows that $G \cdot D = G$. Therefore, $D \in \inv{G}$ by Proposition~\ref{property of xG n > 1} (a).

Next, let $D=(d_{i,j})_{i,j\le n+m} \in \inv{G}$. We need to show that 
$D$ has the form~\eqref{triangulation}, 
where matrices
$A,B,C,O$ satisfy~\eqref{eq:ABCO}.
Since $\inv{G}\subseteq GL(m+n)$, we have
\begin{equation}
\label{one:star}
D \in GL(n+m).
\end{equation}
Since $D\in\inv{G}$, we have
\begin{equation}
\label{two:stars}
G \cdot D = G
\end{equation}
by Proposition~\ref{property of xG n > 1} (a).

Let $\overline{e}_1, \overline{e}_2,\dots,\overline{e}_n$ be the standard basis of $\R^n$; that is, the $i$th coordinate of the vector  $\overline{e}_i$ is equal to $1$ and all its other coordinates are equal to $0$.

Assume that
$1 \leq i \leq n, 1 \leq j \leq n$.
From $\overline{e}_i\in G$ and \eqref{two:stars}, we get
$\overline{e}_i\cdot D\in G\cdot D=G$. 
From this and the definition of $G$, we conclude that the $j$th coordinate $d_{i,j}$ of the vector $\overline{e}_i\cdot D$ must belong to $\Q$.
So  $d_{i,j} \in \mathbb Q$ for all $1 \leq i \leq n, 1 \leq j \leq n$. 

Assume that
$n+1 \leq i \leq n+m, 1 \leq j \leq n$ and $d_{i,j} \ne 0$. Let $y$ be an arbitrary irrational number. 
Then $\overline{x}=y \cdot d_{i,j}^{-1}\cdot \overline{e}_i\in G$, 
so $\overline{x}\cdot D\in G\cdot D=G$ by \eqref{two:stars}. 
From this and the definition of $G$, we conclude that the $j$th coordinate of the vector $\overline{x}\cdot D$ must belong to $\Q$.
On the other hand, this coordinate 
is equal to
$y \notin \mathbb Q$, giving a contradiction.
So $d_{i, j} = 0$ for all $n+1 \leq i \leq n+m, 1 \leq j \leq n$. 

Let $a_{i,j} = d_{i,j}$ when $1 \leq i \leq n, 1 \leq j \leq n$, and   
$c_{i,j} = d_{i+n,j+n}$ when $1 \leq i \leq m, 1 \leq j \leq m$. Then $A = (a_{i,j})_{i,j\le n} \in M_{n,n}$ and $C = (c_{i, j})_{i,j\le m} \in M_{m,m}$. Note that $\dim D \ne 0$  by \eqref{one:star}. Since $\det D = \det A \cdot \det C$, we get $\dim A \ne 0$ and $\dim C \ne 0$. 
This finishes the check of the inclusions $A \in GL_\mathbb Q (n)$ and 
$C \in GL (m)$.
$\Box$ 

\medskip

The group $GL(n)$ has a natural Euclidean topology (coming from $\R^{n^2}$). Therefore, one can discuss the covering dimension $\dim$ of any subgroup of $GL(n)$. (We refer the reader to \cite{E} for a definition of $\dim$ and its properties.)

\begin{theorem}
\label{dimension:theorem}
For every integer $n\ge 1$, define
$D_n = \{\dim \inv{G}: G \mbox{ a dense subgroup of } \mathbb R^n \}$. Then 
\begin{equation}
\label{eq:D_n}
\{qn: q=0,1,\dots,n\}\subseteq D_n\subseteq \{0,1,\dots,n^2\}.
\end{equation}
\end{theorem}
Proof.
Since $\inv{G}\subseteq GL(n)$ and the latter group has a natural Euclidean topology coming from $\R^{n^2}$, it follows that 
$\dim \inv{G}\le \dim \R^{n^2}=n^2$. This establishes the right inclusion in \eqref{eq:D_n}.

Let us check the left inclusion in \eqref{eq:D_n}.

Since $\inv{\Q^n}=GL_{\Q}(n)$ by Corollary~\ref{Aut(Q^n)},
and $GL_{\Q}(n)$ is homeomorphic to an open subset of $\Q^{n^2}$, we have 
$\dim\inv{\Q^n}=0$. Since $\Q^n$ is dense in $\R^n$, this shows that 
$0\in D_n$.

Similarly, $\inv{\R^n}=GL(n)$ by Corollary~\ref{atomorphisms:of:R^n}, and the latter group is homeomorphic to an open subset of $\R^{n^2}$, so $\dim \inv{\R^n}=n^2$.
This proves that $n^2\in D_n$. 

Suppose now that $q$ is an integer satisfying $1\le q\le n-1$.
Then $p=n-q$ is also an integer such that $p\ge 1$. Clearly,
$p+q=n$.
By
Proposition~\ref{Q^n:times:R^m},
$$
\inv{\mathbb Q^p \times \mathbb R^q} = \left\{ \left(\begin{array}{rr}  A & B \\ O & C  \end{array}\right) \in M_{n, n}: A \in GL_\mathbb Q (p), B \in M_{p, q}, C \in GL (q),  O = O_{q, p} \right\}. $$
Since $\inv{\mathbb Q^p \times \mathbb R^q}$ is homeomorphic to an open subset of $\mathbb Q^{p^2}\times \mathbb R^{q n}$, we get $\dim \inv{\mathbb Q^p \times \mathbb R^q} = qn$. Therefore, $qn\in D_n$. 
$\Box$

\medskip
Theorem~\ref{dimension:theorem} motivates the following 

\begin{problem}\label{problem:D_n}
Calculate the set
$D_n$
for every positive integer $n$.
Does the equality $D_n=\{0,1,\dots,n^2\}$ hold?
\end{problem}

In the case $n=2$, it follows from Theorem~\ref{dimension:theorem}
that $\{0,2,4\}\subseteq D_2\subseteq \{0,1,2,3,4\}$.
However, the answer to the following two-dimensional version of 
Problem~\ref{problem:D_n} remains unclear:

\begin{question} Does there exist a dense subgroup $G$ of $\mathbb R^2$ such that $\dim \inv{G} = 1$  or $3$?
\end{question}

It is unclear whether  Corollary~\ref{cardinality:of:automorphism:group:dimension:1} can be extended to the finite-dimensional case.

\begin{question}
\label{Aut:que}
Let $G$ be a dense subgroup of $\R^n$. Is it true that either $G$ is rigid, or the automorphism group $\mathrm{Aut}(G)$ of $G$ is infinite?
\end{question}

As was shown in Corollary~\ref{cardinality:of:automorphism:group:dimension:1}, the answer to 
this question is positive in case $n=1$.

Example~\ref{rigid_subgroups_of_R} motivates the following

\begin{question}
Is there a rigid dense subgroup of $\R^n$ for $n\ge 2$?
\end{question}

\bigskip
\noindent(V.A. Chatyrko)\\
Department of Mathematics, Linkoping University, 581 83 Linkoping, Sweden.\\
vitalij.tjatyrko@liu.se

\vskip0.3cm
\noindent(D.B. Shakhmatov) \\
Graduate School of Science and Engineering, 
Division of Science,
Ehime University\\
Bunkyo-cho 2-5,
790-8577 Matsuyama,
Japan\\
dmitri.shakhmatov@ehime-u.ac.jp
\end{document}